\documentclass{article}
\usepackage{graphicx}

\def\be{\begin{equation}}
\def\ee{\end{equation}}
\def\ff#1{\mbox{\boldmath $#1$} }
\def\a{\alpha}

\def\x{\ff{x}}

\begin{document}

\title{Bat Algorithm is Better Than Intermittent Search Strategy}

\author{Xin-She Yang \\
School of Science and Technology, Middlesex University, \\
London NW4 4BT, United Kingdom.
\and
Suash Deb \\
Department of Computer Science and Engineering, \\
Cambridge Institute of Technology, Ranchi, India. \\
\and
Simon Fong \\
Department of Computer and Information Science, \\
University of Macau, Macau SAR. }

\date{}

\maketitle
\begin{abstract}

The efficiency of any metaheuristic algorithm largely depends on the way of balancing
local intensive exploitation and global diverse exploration. Studies show that
bat algorithm can provide a good balance between these two key components with superior
efficiency. In this paper, we first review some commonly used metaheuristic algorithms,
and then compare the performance of bat algorithm with the so-called
intermittent search strategy. From simulations, we found that  bat algorithm is better than the optimal
intermittent search strategy. We also analyse the comparison results and their implications for
higher dimensional optimization problems. In addition, we also apply bat algorithm
in solving business optimization and engineering design problems.
\end{abstract}

{\bf Citation details:}
X. S. Yang, S. Deb, S. Fong,  Bat Algorithm is Better Than Intermittent Search Strategy, 
{\it Multiple-Valued Logic and Soft Computing}, 22 (3), 223-237 (2014).

\section{Introduction}

Global optimization, computational
intelligence and soft computing often use metaheuristic algorithms.
These algorithms are usually nature-inspired, with multiple interacting agents.
A subset of metaheuristcs are often referred to as swarm intelligence (SI) based
algorithms, and these SI-based algorithms have been developed by mimicking
the so-called swarm intelligence characteristics of biological agents such as
birds, fish, humans and others. For example, particle swarm optimization
was based on the swarming behavior of birds and fish \cite{Kennedy},
while the firefly algorithm was based on the flashing pattern of tropical fireflies \cite{Yang,Yang2}
and cuckoo search algorithm was inspired by the brood parasitism of some cuckoo species \cite{YangDeb}.

In the last two decades, more than a dozen new algorithms such as particle swarm optimization, differential evolution, ant and bee algorithms,
bat algorithm, firefly algorithm and cuckoo search have appeared and they have shown great potential in solving
tough engineering optimization problems \cite{Yang,Floudas,Cui,YangDeb2010,Parp,YangRev2011,YangBook,YangEagle2012,Gandomi,GandomiCAPSO,Cui2}.
In this paper, we provide a first attempt to give some theoretical basis for the
optimal balance of exploitation and exploration for 2D multimodal objective functions.
Then, we use it for choosing algorithm-dependent parameters.
Then, we use bat algorithm as an example to show how such optimal combination can be
achieved.

\section{Brief Review of Metaheuristic Algorithms}

Nature-inspired algorithms, especially those based on swarm intelligence, have
become very popular in recent years. There are a few recent books which are solely
dedicated to metaheuristic algorithms \cite{Yang,Talbi,YangBook,Yangbook2013}.
Metaheuristic algorithms are very diverse, and here we will provide
a brief review of some commonly used metaheuristics.

Ant algorithms, especially ant colony optimization \cite{Dorigo},
mimic the foraging behaviour of social ants, while
bees-inspired algorithms use various characteristics
such as waggle dance, polarization and nectar maximization
are often used to simulate the allocation of the foraging bees along flower patches
and thus different search regions in the search space,  depending on many factors such as
the nectar richness and the proximity to the bee hive \cite{Nak,Yang05,Kara,Pham,Afshar}.
For a more comprehensive review,
please refer to Yang \cite{YangBook} and Parpinelli and Lope \cite{Parp}.

Genetic algorithms are a class of algorithms using genetic operators such as crossover, mutation and selection
of the fittest,  pioneered by J. Holland and his collaborators in the 1960s
and 1970s \cite{Holland}. On the other hand, differential evolution (DE) was developed by R. Storn and K. Price by their
nominal papers  in 1996 and 1997 \cite{Storn,StornPrice,Price}.

Particle swarm optimization (PSO) was developed by Kennedy and
Eberhart in 1995 \cite{Kennedy}, based on the swarm behaviour such
as fish and bird schooling in nature. There are many variants which extend the standard PSO
algorithm \cite{Kennedy,Yang,YangBook,GandomiCAPSO}.

Firefly Algorithm (FA) was first developed by Xin-She Yang in 2007 \cite{Yang,Yang2,YangFA,GandomiCFA}
which was based on the flashing patterns and behaviour of fireflies.
Firefly algorithm has attracted much attention \cite{Sayadi,Apo,GandomiFA}.
A discrete version of FA can efficiently solve NP-hard scheduling problems \cite{Sayadi},
while a detailed analysis has demonstrated the efficiency of FA over a wide range of test problems,
including multobjective load dispatch problems \cite{Apo}.
High nonlinear and non-convex global optimization problems can be solved
using firefly algorithm efficiently \cite{GandomiFA}.

Cuckoo search (CS) is one of the latest nature-inspired metaheuristic algorithms, developed in 2009
by Xin-She Yang and Suash Deb \cite{YangDeb}.
CS is based on the brood parasitism of some cuckoo species. Furthermore, CS can be further enhanced
by L\'evy flights \cite{Pav}.
Extensive studies show that CS can perform much better than PSO and genetic algorithms.
There have been a lot of attention and recent studies using cuckoo search  with
diverse range of applications \cite{Gandomi,Walton,Durgun}.
Walton et al. improved the algorithm by formulating a modified cuckoo search algorithm \cite{Walton},
while Yang and Deb extended it to multiobjective optimization problems
by Yang and Deb in 2011, though the paper was published later in 2013 \cite{YangDeb2011}.

\section{Balancing Exploration and Exploitation}

The main components of any metaheuristic algorithms are: intensification
and diversification, or exploitation and exploration \cite{Blum,Yang12}.
Diversification can generate diverse, explorative solutions, and intensification
provides a mechanism for exploiting local information so as to speed up
convergence. In addition, selection of the best solutions often
ensure the quality solution will remain in the population so as to avoid potential divergence.

Exploration in metaheuristics can be achieved often by
the use of randomization \cite{Blum,Yang,YangFA},
which enables an algorithm to have the ability to jump out of any local optimum so as to explore the
search space globally. Randomization often use random walks, and it can
use for both local and global search, depending on the step sizes and
way to generate new solutions. Fine-tuning the right amount of randomness
and balancing local search and global search may be essential for any
algorithm to perform well.

Exploitation is the use of local knowledge of the search and solutions found so far
so that new search moves can concentrate on the local regions or neighborhood
where the optimality may be close; however, this local optimum may not be the
global optimality. Exploitation tends to use strong local information such as gradients,
the shape of the mode such as convexity, and the history of the search process.
A classic technique is the so-called hill-climbing which uses the local gradients
or derivatives intensively.

Empirical knowledge from observations and simulations of the convergence behaviour of
common optimization algorithms suggests that exploitation tends to increase the
speed of convergence, while exploration tends to decrease the convergence rate of the
algorithm. On the other hand, too much exploration increases the probability of finding the global
optimality, while strong exploitation tends to make the algorithm being trapped in a
local optimum. Therefore, there is a fine balance between the right amount of exploration
and the right degree of exploitation. Despite its importance, there is no practical guideline
for this balance.

\section{Intermittent Search Strategy}

Even there is no guideline in practice, some preliminary work on the very limited cases exists in the
literature and may provide some insight into the possible choice of parameters so as to balance these components.
Intermittent search strategies concern an iterative strategy consisting of a slow phase and a fast phase \cite{Ben,Ben2}.
Here the slow phase is the detection phase by slowing down and intensive, static local search techniques,
while the fast phase is the search without detection and can be considered as an exploration technique.
For example, the static target detection with a small region of radius $a$ in a much larger region $b$
where $a \ll b$ can be modelled as a slow diffusive process in terms of random walks with a diffusion
coefficient $D$.

Let $\tau_a$ and $\tau_b$ be
the mean times spent in intensive detection stage and the time spent in the exploration stage,
respectively, in the 2D case. The diffusive search process is governed by
the mean first-passage time satisfying the following equations
\be D \nabla_r^2 t_1 + \frac{1}{2 \pi \tau_a} \int_0^{2 \pi} [t_2(r) - t_1(r)] d\theta+1=0, \ee
\be \ff{u} \cdot \nabla_r t_2(r) -\frac{1}{\tau_b} [t_2(r)-t_1(r)]+1=0, \ee
where $t_2$ and $t_1$ are mean first-passage times during the
search process, starting from slow and fast stages, respectively, and $\ff{u}$ is the mean search speed \cite{Ben2}.

After some lengthy mathematical analysis, the optimal balance of these two stages can be estimated as
\be r_{\rm optimal}=\frac{\tau_a}{\tau_b^2} \approx \frac{D}{a^2} \frac{1}{[2-\frac{1}{\ln(b/a)}]^2}. \label{equ-balance} \ee
Assuming that the search steps have a uniform velocity $u$ at each step on average, the minimum times required
for each phase can be estimated as
\be \tau_a^{\min} \approx \frac{D}{2 u^2} \frac{\ln^2(b/a)}{[2 \ln(b/a)-1]}, \ee
and
\be \tau_b^{\min} \approx \frac{a}{u} \sqrt{\ln(b/a)-\frac{1}{2}}. \ee
When $u \rightarrow \infty$, these relationships lead to the above optimal ratio of two stages.

It is worth pointing out that the above result is only valid for 2D cases, and there is no general results
for higher dimensions, except in some special 3D cases \cite{Ben}. Now let us use this limited results to help
choose the possible values of algorithm-dependent parameters in bat algorithm \cite{YangBA}, as an example.

For higher-dimensional problems, no result exists. One possible extension is to
use extrapolation to get an estimate. Based on the results on 2D and 3D cases \cite{Ben2},
we can estimate that for any $d$-dimensional cases $d \ge 3$
\be \frac{\tau_1}{\tau_2^2} \sim O\Big(\frac{D}{a^2}\Big), \quad \tau_m \sim O\Big( \frac{b}{u} (\frac{b}{a})^{d-1} \Big), \ee
where $\tau_m$ the  mean search time or average number of iterations.
This extension may not be good news for higher dimensional problems, as the mean number of function evaluations to
find optimal solutions can increase exponentially as the dimensions increase. However, in practice,
we do not need to find the guaranteed global optimality, we may be satisfied with suboptimality,
and sometimes we may be lucky to find such global optimality even with a limited/fixed number of iterations.
This may indicate there is a huge gap between theoretical understanding and the observations as well as run-time
behaviour in practice. More studies are highly needed to address these important issues.

\section{Bat Algorithm}

Bat algorithm was developed by Xin-She Yang in 2010\cite{YangBA},
which were based on the fascinating characteristics of echolocation of microbats \cite{Bat}.
Bat algorithm uses the following two idealized rules:

\begin{itemize}

\item Bats fly randomly with velocity $\ff{v}_i$
at position $\x_i$ to search for food/prey. They emit short pulses at a rate of $r$,
with a varying frequency $f_{\min}$, wavelength $\lambda$
and loudness $A_0$.

\item Frequency (or wavelength) and the rate of pulse emission $r \in [0,1]$ can be
adjusted,  depending on the proximity of their target.

\end{itemize}

Loudness can vary in many ways, though a common scheme is to vary from
a large (positive) $A_0$ to a minimum constant value $A_{\min}$.
For simplicity, no ray tracing is used. Furthermore,  the frequency $f$  in
a range $[f_{\min}, f_{\max}]$ corresponds to a range
of wavelengths $[\lambda_{\min}, \lambda_{\max}]$.

\subsection{Bat Motion}

The  positions $\x_i$ and velocities $\ff{v}_i$ of bat $i$
in a $d$-dimensional search space are updated as follows:
\be f_i =f_{\min} + (f_{\max}-f_{\min}) \beta,\ee
\be \ff{v}_i^{t+1} = \ff{v}_i^{t} +  (\x_i^t - \x_*) f_i, \ee
\be \x_i^{t+1}=\x_i^{t} + \ff{v}_i^t,  \ee
where the random number is drawn uniformly from $\beta \in [0,1]$.
Here $\x_*$ is the current global best
location (solution) which is located after comparing all
the solutions among all the $n$ bats at each iteration $t$.
For ease of implementation, $f_{\min}=0$ and $f_{\max}=O(1)$ are used,
depending on the domain size of the problem of interest.
At the start of the simulation, each bat is randomly initialized with
a frequency in $[f_{\min}, f_{\max}]$.

As part of local random walks, once a solution is selected among the
current best solutions, a new solution for each bat is generated by
\be \x_{\rm new}=\x_{\rm old} + \ff{\epsilon}  \; A^{t}, \ee
where $\ff{\epsilon}$ is a random number vector drawn from $[-1,1]$, while $A^{t}=<\!\!A_i^{t}\!>$
is the average loudness of all the bats at this time step.

Furthermore, for simplicity, we can also use
$A_0=1$ and $A_{\min}=0$, assuming $A_{\min}=0$ means that a bat
has just found the prey and temporarily stop emitting any sound.
Now we have \cite{YangBA}
\be A_i^{t+1}=\alpha A_{i}^{t}, \;\;\;\;\; r_i^{t}= r_i^0 [1-\exp(-\gamma t)],
\label{rate-equ-50} \ee
where $\alpha$ and $\gamma$ are constants.
Here, $\alpha$ has a similar role to that played by the
cooling factor of a cooling schedule in the simulated annealing.
For any $0<\alpha<1$ and $\gamma>0$, we have
\be A_i^t \rightarrow 0, \;\;\; r_i^t \rightarrow r_i^0, \;\;\textrm{as} \;\;
t \rightarrow \infty. \ee
Obviously,  $\a=\gamma$ can be assumed for simplicity, and
in fact $\a=\gamma=0.9$ has been used in our simulations.
Bat algorithm has been proved to be very efficient for solving various optimization problems \cite{YangGandomiBA}.

\section{Numerical Experiments}

\subsection{Landscape-Dependent Optimality}

If we use the 2D simple, isotropic random walks for local exploration, then we have
\be D\approx \frac{s^2}{2}, \ee
where $s$ is the step length with a jump during a unit time interval or each iteration step.
From equation (\ref{equ-balance}), the optimal ratio of exploitation and exploration
in a special case of $b \approx 10 a$ becomes
\be \frac{\tau_a}{\tau_b^2} \approx 0.2. \ee
In case of $b/a \rightarrow \infty$, we have $\tau_a/\tau_b^2 \approx 1/8.$
which implies that more times should spend on the exploration stage. It is worth pointing out that
the naive guess of 50-50 probability in each stage is not the best choice. More efforts
should focus on the exploration so that the best solutions found by the algorithm can be
globally optimal with possibly the least computing effort.

In the case studies to be described below, we have used the bat algorithm to find the
optimal solutions to two design benchmarks. If we set $\tau_b=1$ as the reference timescale, then
we found that the optimal ratio is between 0.18 to 0.24, which are roughly close to the above theoretical result.
This may imply that bat algorithm has an intrinsic ability of balancing exploration
and exploitation close to the true optimality.

\subsection{Standing-Wave Function}

Let us first use a multimodal test function to see how to find the fine balance between
exploration and exploitation in an algorithm for a given task. A standing-wave test function
can be a good example \cite{YangFA}
\be f(\x)=1+ \Big\{\exp[-\sum_{i=1}^d (\frac{x_i}{\beta})^{10}]  - 2 \exp[-\sum_{i=1}^d (x_i-\pi)^2]  \Big\} \cdot \prod_{i=1}^d \cos^2 x_i, \ee
which is multimodal with many local peaks and valleys. It has a unique global minimum at
$f_{\min}=0$ at $(\pi,\pi,...,\pi)$ in the domain $-20 \le x_i \le 20$ where $i=1,2,...,d$ and $\beta=15$.
In this case, we can estimate that $R=20$ and $a \approx \pi/2$, this means that $R/a \approx 12.7$, and
we have in the case of $d=2$
\be p_e \approx \tau_{\rm optimal} \approx \frac{1}{2 [2-1/\ln (R/a)]^2} \approx 0.19. \ee
This indicate that the algorithm should spend 80\% of its computational effort on global explorative search,
and 20\% of its effort on local intensive search.

For the bat algorithm, we have used $n=15$ and $1000$ iterations.
We have calculated the fraction of iterations/function evaluations
for exploitation to exploration, that is $Q=$ exploitation/exploration,
thus $Q$ can thus affect the solution quality. A set of 25 numerical experiments
have been carried out for each value of $Q$ and the results are summarized in Table 1.

\begin{table}
\begin{center}
\caption{Variations of $Q$ and its effect on the solution quality. }
\begin{tabular}{|l|l|l|l|l|l|l|l|l}
\hline
$Q \!\!$  & 0.3 & 0.2  & 0.1 & $0.05 \!\!$ \\ \hline
$f_{\min}\!\!$  & 1.1e-12 & 2.3e-14 & 7.7e-12 & 8.1e-11 \\
\hline
\end{tabular}
\end{center}
\end{table}

This table clearly shows that $Q \approx 0.2$ provides the optimal balance of local exploitation
and global exploration, which is consistent with the theoretical estimation.

Though there is no direct analytical results for higher dimensions, we can expect that
more emphasis on global exploration is also true for higher dimensional optimization problems.
Let us study this test function for various higher dimensions.

\subsection{Comparison for Higher Dimensions}

As the dimensions increase, we usually expect the number of iterations of finding the global
optimality should increase. In terms of mean search time/iterations,
B\'enichou et al.'s intermittent search theory suggests that
\be \tau_m\Big|_{(d=1)} = \frac{2b}{u} \sqrt{\frac{b}{3a}}, \ee
\be \tau_m \Big|_{(d=2)} =\frac{2 b^2}{a u} \sqrt{\ln (b/a)}, \ee
\be \tau_m\Big|_{(d=3)}=\frac{2.2 b}{u} (\frac{b}{a})^2. \ee

For higher dimensions, we can only estimate the main trend based on the intermittent search strategy. That is,
\be \frac{\tau_1}{\tau_2^2} \sim O\Big(\frac{D}{a^2}\Big), \quad \tau_m \sim O\Big( \frac{b}{u} (\frac{b}{a})^{d-1} \Big), \ee
which means that number of iterations may increase exponentially with the dimension $d$.
It is worth pointing out that the ratio between the two stage are independent of the dimensions.
In other words, once  we find the optimal balance between exploration and exploitation, we can use the
algorithm for any high dimensions. Now let us use bat algorithm to carry out search in higher dimensions
for the above standing wave function and compare its performance with the implication of intermittent search strategy.

\begin{figure}
\centerline{\includegraphics[height=2.5in,width=3in]{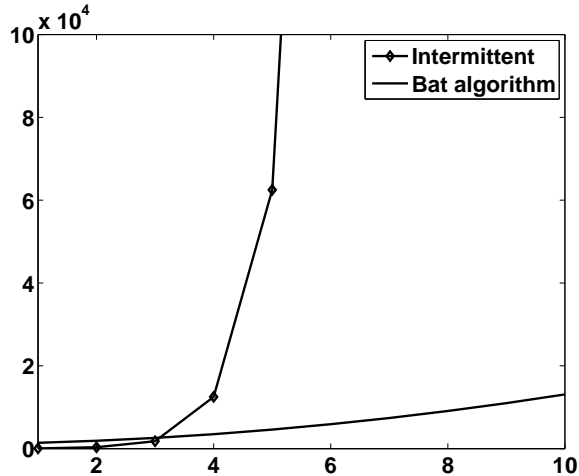}}
\caption{Comparison of the actual number of iterations with the theoretical results by the intermittent search strategy.
This clearly shows that bat algorithm is better than the intermittent search strategy. \label{fig-bat}}
\end{figure}

For the case of $b=20$, $a=\pi/2$ and $u=1$, Fig. \ref{fig-bat} shows the comparison of the numbers of iterations suggested by
intermittent search strategy and the actual numbers of iterations using bat algorithm to obtain the
globally optimal solution with a tolerance or accuracy of 5 decimal places.
It can be seen clearly that the number of iterations needed by the intermittent search strategy increases
exponentially versus the number of dimensions, while the actual number of iterations used in the algorithm
only increases slightly, seemingly weakly a low-order polynomial. This suggests that bat algorithm is very efficient
and requires far fewer (and often many orders lower) number of iterations.

\section{Business Optimization and Engineering Applications}

\subsection{Project Scheduling}

In order to validate and test the proposed method, we use the resource-constrained project scheduling problems
by Kolisch and Sprecher \cite{Kol,Kol2}. The basic model consists
of $J$ activities/tasks, and some activities cannot start before all its predecessors $h$
are completed. In addition, each activity $j=1,2,...,J$
can be carried out, without interruption,
in one of the $M_j$ modes, and performing any activity $j$ in any chosen
mode $m$ takes $d_{jm}$ periods, which is supported by a set of renewable
resource $R$ and non-renewable resources $N$. The project's makespan or upper
bound is T, and the overall capacity of non-renewable resources is
$K_r^{\nu}$ where $r \in N$. For an activity $j$ scheduled in mode $m$,
it uses $k^{\rho}_{jmr}$ units of renewable resources
and $k^{\nu}_{jmr}$ units of non-renewable resources
in period $t=1,2,..., T$.

\begin{table}[ht]
\caption{Kernel parameters used in SVM.}
\centering
\begin{tabular}{lllll}
\hline \hline
Iterations & Kernel parameters \\
\hline
1000 & $C=149.2$, $\sigma^2=67.9$ \\
5000 & $C=127.9$, $\sigma^2=64.0$ \\
\hline
\end{tabular}
\end{table}

Using the online benchmark library \cite{Kol2}, we have solved this type
of problem with $J=30$ activities (the standard test set j30).  The run time
on a modern desktop computer is about 3.1 seconds for $N=1000$ iterations
to 17.2 seconds for $N=5000$ iterations. We have run
the simulations for 50 times so as to obtain meaningful statistics.
The main results  are analyzed in Table 3.

\begin{table}[ht]
\caption{Mean deviations from the optimal solution (J=30).}
\centering
\begin{tabular}{lllll}
\hline \hline
Algorithm  & $N=1000$ & $5000$ \\
\hline
PSO  & 0.26 & 0.21 \\
hybribd GA  & 0.27 & 0.06 \\
Tabu search   & 0.46 & 0.16 \\
Adapting GA   & 0.38 & 0.22 \\
{\it Current BA} & {\it 0.29} & {\it 0.049} \\
\hline
\end{tabular}
\end{table}

\subsection{Spring Design}

In engineering optimization, a well-known design benchmark is
the spring design problem with three design variables: the wire diameter $w$,
the mean coil diameter $d$, and the length (or number of coils) $L$.
For detailed description, please refer to earlier
studies \cite{Gandomi,Gandomi2}.
This problem can be written compactly as
\be \textrm{Minimise } f(\x)=(L+2) w^2 d, \qquad \qquad  \ee
subject to
\be \begin{array}{lll}
g_1(\x) =1-\frac{d^3 L}{71785 w^4} \le 0, \\ \\
g_2(\x) =1-\frac{140.45 w}{d^2 L} \le 0, \\ \\
g_3(\x) =\frac{2(w + d)}{3}-1 \le 0, \\ \\
g_4(\x)= \frac{d (4 d-w)}{w^3 (12566  d - w)} + \frac{1}{5108 w^2} -1 \le 0,
\end{array} \ee
with the following limits
\be 0.05 \le w \le 2.0, \quad 0.25 \le d \le 1.3, \quad 2.0 \le L \le 15.0. \ee

Using bat algorithm, we have obtained
\be f_*=0.012665  \quad \textrm{ at } \;\; (0.051690, 0.356750, 11.287126). \ee

\subsection{Welded Beam Design}

The welded beam design is another well-known benchmark \cite{Yang,Gandomi,Gandomi2}.
The problem has four design variables: the width $w$
and length $L$ of the welded area, the depth $h$ and thickness $h$ of the main
beam. The objective is to minimise the overall fabrication cost, under
the appropriate constraints of shear stress $\tau$,
bending stress $\sigma$, buckling load $P$ and maximum end deflection $\delta$.

The problem can be written as
\be \textrm{minimise } \; f(\x)=1.10471 w^2 L + 0.04811 d h (14.0+L), \ee
subject to
\be
\begin{array}{lll}
g_1(\x)=w -h \le 0, \vspace{2pt} \\ \vspace{3pt}
g_2(\x) =\delta(\x) - 0.25 \le 0, \\ \vspace{3pt}
g_3(\x)=\tau(\x)-13,600 \le 0, \\ \vspace{3pt}
g_4(\x)=\sigma(\x)-30,000 \le 0, \\ \vspace{3pt}
g_5(\x)=0.10471 w^2 +0.04811  h d (14+L) -5.0 \le 0, \\ \vspace{3pt}
g_6(\x)=0.125 - w \le 0, \\ \vspace{3pt}
g_7(\x)=6000 - P(\x) \le 0,
\end{array}
\ee
where
\be \begin{array}{ll}
 \sigma(\x)=\frac{504,000}{h d^2},  & Q=6000 (14+\frac{L}{2}), \\ \\
 D=\frac{1}{2} \sqrt{L^2 + (w+d)^2}, & J=\sqrt{2} \; w L [ \frac{L^2}{6} + \frac{(w+d)^2}{2}], \\ \\
 \delta=\frac{65,856}{30,000 h d^3}, &  \beta=\frac{QD}{J}, \\ \\
 \alpha=\frac{6000}{\sqrt{2} w L}, & \tau(\x)=\sqrt{\alpha^2 + \frac{\alpha \beta L}{D}+\beta^2}, \\ \\
 P=0.61423 \times 10^6 \; \frac{d h^3}{6} (1-\frac{d \sqrt{30/48}}{28}). &
\end{array} \ee
The simple bounds for design variables are $0.1 \le L, d \le 10$ and
$0.1 \le w, h \le 2.0$.

Using our bat algorithm, we have obtained
\[ \x_*=(w ,L, d, h) \] \be =
(0.20572963978,   3.47048866563,   9.03662391036,   0.20572963979), \ee
which gives the optimal solution \be f(\x*)_{\min} = 1.724852308598.      \ee

\section{Conclusions}

Nature-inspired metaheuristic algorithms have gained popularity, which is partly due to
its ability of dealing with nonlinear global optimization problems. We have
highlighted the importance of exploitation and exploration and their effect on the
efficiency of an algorithm. Then, we use the intermittent search strategy theory as a preliminary
basis for analyzing these key components and ways to find the possibly optimal settings
for algorithm-dependent parameters.

In addition, we have used the bat algorithm
to find this optimal balance, and confirmed that bat algorithm can indeed
provide a good balance of exploitation and exploration.
We have also shown that
bat algorithm requires far fewer function evaluations. However, the huge differences between
intermittent search theory and the behaviour of metaheuristics in practice also suggest
there is still a huge gap between our understanding of algorithms and the actual behaviour of
metaheuristics. It is highly desirable to carry out more research in this area.

\end{document}